%% file: smooth-bounded-sums-main.tex
\documentclass[12pt]{article}
\usepackage[nottoc]{tocbibind}
\usepackage[latin1]{inputenc}
\usepackage{fancyhdr}
\usepackage{verbatim}
\usepackage{indentfirst}
\usepackage{graphicx}
\usepackage{epstopdf}
\usepackage{setspace}
\usepackage{color}
\usepackage{newlfont}
\usepackage{amssymb}
\usepackage[intlimits]{amsmath}
\usepackage{latexsym}
\usepackage{amsthm}
\usepackage{amscd}
\usepackage{mathrsfs}
\usepackage[dvips, outer=1in, top=1in, bottom=1in]{geometry}
\usepackage{mathtools, array, multirow}
\usepackage{hyperref}
\usepackage{rotating}
\usepackage{mathtools}
\usepackage{dsfont}

\theoremstyle{plain}                    
\newtheorem{theorem}{Theorem}[section]
\newtheorem*{theorem*}{Theorem}
\newtheorem{lemma}[theorem]{Lemma}    
\newtheorem{prop}[theorem]{Proposition}
\newtheorem{proposition}[theorem]{Proposition}
\newtheorem{cor}[theorem]{Corollary}
\theoremstyle{definition}

\newtheorem{ex}[theorem]{Example}        
\newtheorem{remark}[theorem]{Remark} 

\newcommand{\N}{\mathbb{N}}
\newcommand{\Z}{\mathbb{Z}}
\newcommand{\Q}{\mathbb{Q}}
\newcommand{\R}{\mathbb{R}}
\newcommand{\C}{\mathbb{C}}

\newcommand{\T}{\mathbb{T}}
\newcommand{\bianco}{\textcolor{white}{.}}
\newcommand{\uhp}{\mathfrak{H}}

\newcommand{\h}{\mathfrak{H}}
\newcommand{\sltr}{\mathrm{SL}(2,\R)}
\newcommand{\sltz}{\mathrm{SL}(2,\Z)}
\newcommand{\asltr}{\mathrm{ASL}(2,\R)}

\newcommand{\tsltr}{\widetilde{\mathrm{SL}}(2,\R)}

\newcommand{\ha}{\frac{1}{2}}
\newcommand{\tha}{\tfrac{1}{2}}

\newcommand{\de}{\mathrm{d}}
\newcommand{\psltr}{\mathrm{PSL}(2,\R)}

\newcommand{\Hei}{\mathbb{H}(\mathbb{R})}

\newcommand{\Heik}{\mathbb{H}(\mathbb{R}^k)}
\newcommand{\ba}{\begin{array}}
\newcommand{\ea}{\end{array}}
\newcommand{\ve}[2]{\left(\ba{c}\!#1\!\\ \!#2\!\ea\right)}
\newcommand{\sve}[2]{\left(\begin{smallmatrix}\!#1\!\\ \!#2\!\end{smallmatrix}\right)}
\newcommand{\veH}[3]{\left(\!\begin{pmatrix}#1\\ #2\end{pmatrix},#3\right)}
\newcommand{\sveH}[3]{\left(\left(\begin{smallmatrix}\!#1\!\\ \!#2\!\end{smallmatrix}\right),#3\right)}

\newcommand{\bo}[1]{\mathbf{#1}}

\newcommand{\LtRk}{\mathrm L^2(\R^k)}

\newcommand{\sma}[4]{\left(\begin{smallmatrix} #1&#2\\#3&#4\end{smallmatrix}\right)}
\newcommand{\e}[1]{\mathrm{e}\!\left(#1\right)}

\newcommand{\Leb}{\mathrm{Leb}}

\def\veca{{\text{\boldmath$a$}}}
\def\vecb{{\text{\boldmath$b$}}}

\def\vecm{{\text{\boldmath$m$}}}
\def\vecn{{\text{\boldmath$n$}}}

\def\vecs{{\text{\boldmath$s$}}}

\def\vecu{{\text{\boldmath$u$}}}
\def\vecv{{\text{\boldmath$v$}}}

\def\vecw{{\text{\boldmath$w$}}}

\def\vecz{{\text{\boldmath$z$}}}

\def\vecalf{{\text{\boldmath$\alpha$}}}
\def\vecbeta{{\text{\boldmath$\beta$}}}

\def\vecxi{{\text{\boldmath$\xi$}}}

\def\Re{\operatorname{Re}}
\def\Im{\operatorname{Im}}

\def\ASL{\ASL(2,\mathbb{R})}

\def\Onder#1#2#3#4#5{#1 \setbox0=\hbox{$#1$}\setbox1=\hbox{$#2$}
       \dimen0=.5\wd0 \dimen1=\dimen0 \dimen2=\dp0 \dimen3=\dimen2
       \advance\dimen0 by .5\wd1 \advance\dimen0 by -#4
       \advance\dimen1 by -.5\wd1 \advance\dimen1 by -#4
       \advance\dimen2 by -#3 \advance\dimen2 by \ht1
       \advance\dimen2 by 0.3ex \advance\dimen3 by #5
        \kern-\dimen0\raisebox{-\dimen2}[0ex][\dimen3]{\box1}
       \kern\dimen1}
\newcommand{\GaG}{\Gamma\backslash G}

\newcommand{\DaG}{\Delta \backslash G}

\newcommand{\Si}{\mathcal{S}}
\newcommand{\bn}{\mathbf{0}}

\newcommand{\tG}{\widetilde{G}}

\newcommand{\ttM}{\widetilde{M}}

\newcommand{\vtheta}{\vartheta}

\newcommand{\F}{\mathcal{F}}

\newcommand{\matr}[4]{\left( \begin{matrix} #1 & #2 \\ #3 & #4 \end{matrix} \right) }
\newcommand{\smatr}[4]{\bigr( \begin{smallmatrix} #1 & #2 \\ #3 & #4 \end{smallmatrix} \bigr) }

\newcommand{\calC}{\mathcal{C}}

\newcommand{\calA}{\mathcal{A}}
\newcommand{\calK}{\mathcal{K}}

\newcommand{\Thetapair}[2]{\Theta_{#1}\overline{\Theta_{#2}}}
\newcommand{\ind}{\mathds{1}}

\newcommand{\Orb}{\mathrm{Orb}}

\newcommand{\tR}{\widetilde{R}}
\newcommand{\calF}{\mathcal{F}}

\newcommand{\Fone}{\mathcal{F}_{\Gamma}^{(1)}}
\newcommand{\Finfty}{\mathcal{F}_{\Gamma}^{(\infty)}}

\newcommand{\vepsilon}{\varepsilon}

\newcommand{\toab}{^{(\alpha,\beta)}}
\newcommand{\tovab}{^{(\vecalf, \vecbeta)}}

\newcommand{\lp}{\left(}
\newcommand{\rp}{\right)}
\newcommand{\labs}{\left|}
\newcommand{\rabs}{\right|}
\newcommand{\ls}{\left[}
\newcommand{\rs}{\right]}
\newcommand{\lcur}{\left\{}
\newcommand{\rcur}{\right\}}

\numberwithin{equation}{section}

\input xy
\xyoption{all}
\usepackage{import}

\newcommand{\footremember}[2]{%
    \footnote{#2}
    \newcounter{#1}
    \setcounter{#1}{\value{footnote}}%
}

\title{Bounds for Smooth Theta Sums\\ with Rational Parameters}
\author{Francesco Cellarosi\footremember{queens}{Department of Mathematics and Statistics. Queen's University. Kingston, ON, Canada.} and Tariq Osman\footremember{brandeis}{Department of Mathematics. Brandeis University. Waltham, MA, U.S.A.}.}
\date{\today}

\begin{document}

\maketitle

\begin{abstract}
     We provide an explicit family of pairs $(\vecalf, \vecbeta) \in \R^k \times \R^k$ such that for sufficiently regular $f$, there is a constant $C$ for which the theta sum bound
    \begin{align*}
        \labs\sum_{\vecn\in \Z^k} f\!\left(\tfrac{1}{N}\vecn\right)\exp\{2 \pi i((\tha\|\vecn\|^2 + \vecbeta\cdot \vecn)x + \vecalf\cdot\vecn)\}\rabs \leq C N^{k/2},
    \end{align*}
     holds for every $x \in \R$ and every $N \in \N$. Central to the proof is realising that, for fixed $N$, the theta sum  normalised by $N^{k/2}$ agrees with an automorphic function $|\Theta_f|$ evaluated along a special curve known as a horocycle lift. The lift depends on the pair $(\vecalf,\vecbeta)$, and so the bound follows from showing that there are pairs such that $|\Theta_f|$ remains bounded along the entire horocycle lift. 
\end{abstract}

\tableofcontents
\import{./}{introduction.tex}
\import{./}{background_material-01.tex}

\import{./}{uniform-bounds-v1.tex}

\import{./}{key-ingredients.tex}

\import{./}{study-of-the-orbits.tex}

\import{./}{theta-function-bounds.tex}

\bibliographystyle{plain}
\bibliography{thesis-bibliography}

\end{document}

%% file: introduction.tex
\section{Introduction}\label{introduction}
Let $\e{z} := e^{2\pi i z}$. Let $f: \R^k \to \R$ be a Schwartz function. We consider the \emph{generalised quadratic Weyl sum} (or \emph{generalised theta sums}) 
\begin{align}\label{def:S_N^f(x,alpha,beta)-intro}
    S^f_N(x; \vecalf, \vecbeta) := \sum_{\vecn\in \Z^k} f\!\left(\tfrac{1}{N}\vecn\right)\e{(\tha\|\vecn\|^2 + \vecbeta\cdot \vecn)x + \vecalf\cdot\vecn},
\end{align}
where $x\in\R$, $N\in\N$, and $\vecalf, \vecbeta \in \R^k$. We can think of $f$ as a smooth cut-off function.
We prove the following theorem.

\begin{theorem}[Main Theorem, for Schwartz cut-offs.]\label{thm:main-thm-intro}
    Let $f \in \Si (\R^k)$ be a Schwartz function and let $\vecalf, \vecbeta \in \R^k$. If  $\vecalf = (\alpha_1, \cdots, \alpha_n)$,  $\vecbeta = (\beta_1, \cdots, \beta_n)$, and at least one of $(\alpha_i,\beta_i) = (\tfrac{a}{2m}, \tfrac{b}{2m}) \in \Q^2$ is such that $\gcd(a,b, m) = 1$, and $a$, $b$, and $m$ are all odd, then
    \begin{align}\label{eq:statement-main-thm-intro}
        \left|S_N^f(x; \vecalf, \vecbeta) \right| \ll_{k,m,\vecbeta,f}N^{k/2},
    \end{align}
    for every $x \in \R$ and every $N \in \N$. 
\end{theorem}
In \eqref{eq:statement-main-thm-intro} and in the rest of the paper we use Vinogradov's ``$\ll$'' notation (which is equivalent to Landau's $O$-notation) and stress the dependence of the implied constants upon the parameters written as subscripts.
As we shall see, we can relax the assumption that $f$ is Schwartz, see Theorem \ref{thm:main-theorem}. 
When $k = 1$ we have the following corollary. 
\begin{cor}\label{cor:main-cor-intro}
    Let $f \in \Si(\R)$. Let $\alpha = \tfrac{a}{2m}$ and $\beta = \tfrac{b}{2m}$, with $a$, $b$, and $m$ all odd, and such that $\gcd(a,b,m) = 1$. Then
    \begin{align}
        |S_N^f(x; \alpha, \beta)| \ll_{m,\beta,f} \sqrt{N},\label{cor:main-cor-intro-bound}
    \end{align}
      for every $x \in \R$ and every $N \in \N$.
\end{cor}

\begin{remark}
    The only previously known instance of Corollary \ref{cor:main-cor-intro} is due to Marklof when $m=1$, i.e. $(\alpha,\beta)=(\frac{1}{2},\frac{1}{2})$, see Section 5.2 of \cite{Marklof2007b}. Furthermore, it follows from Theorem 1.4 (i) in \cite{Cellarosi-Marklof} and Theorem 1.0.7 in \cite{Cellarosi-Osman-rational-tails} that the limiting distribution of $|N^{-\tha}S_N (x; \alpha,\beta)|$ as $N \to \infty$ for any pair $(\alpha, \beta) \in \R^2$ \emph{not} of the form given in Corollary \ref{cor:main-cor-intro}, is heavy tailed. Therefore, the rational pairs $(\alpha,\beta)$ in  Corollary \ref{cor:main-cor-intro} are \emph{the only pairs} in $\R^2$ for which the bound \eqref{cor:main-cor-intro-bound} holds for every $x \in \R$ and every $N \in \N$. 
\end{remark}

In Section \ref{jacobi-theta-bound-section}, as an illustration of the main theorem, we also obtain bounds for the classical Jacobi theta function
\begin{align}\label{def-classical-jacobi-theta3}
    \vartheta_3(z;w) = \sum_{n \in \Z} e^{2 \pi i nz +  \pi n^2 w},
\end{align}
where $z \in \C$ and $w \in \uhp := \{ x + iy \in \C : y > 0\}$. When $w$ approaches the boundary in $\uhp$ and $z$ is of the form $\alpha+\beta\Re w$ with $(\alpha,\beta)$ as in Corollary \ref{cor:main-cor-intro}, we have the following 

\begin{theorem}\label{thm:jacobi-theta-bound-intro}
    Suppose $(\alpha,\beta) = (\tfrac{a}{2m},\tfrac{b}{2m})\in\R^2/\Z^2$ with $\gcd(a,b,m) = 1$ and $a, b$ and $m$ all odd. 
    Then 
    \begin{align}
    \labs\vtheta_3 (\alpha + \beta x, x + i \vepsilon)\rabs \ll_{m} \vepsilon^{-1/4}.
    \label{statement-thm-bound-theta3}
    \end{align}
\end{theorem}

\begin{figure}[ht!]
    \centering
\hspace{-.6cm}\includegraphics[width=16.35cm]{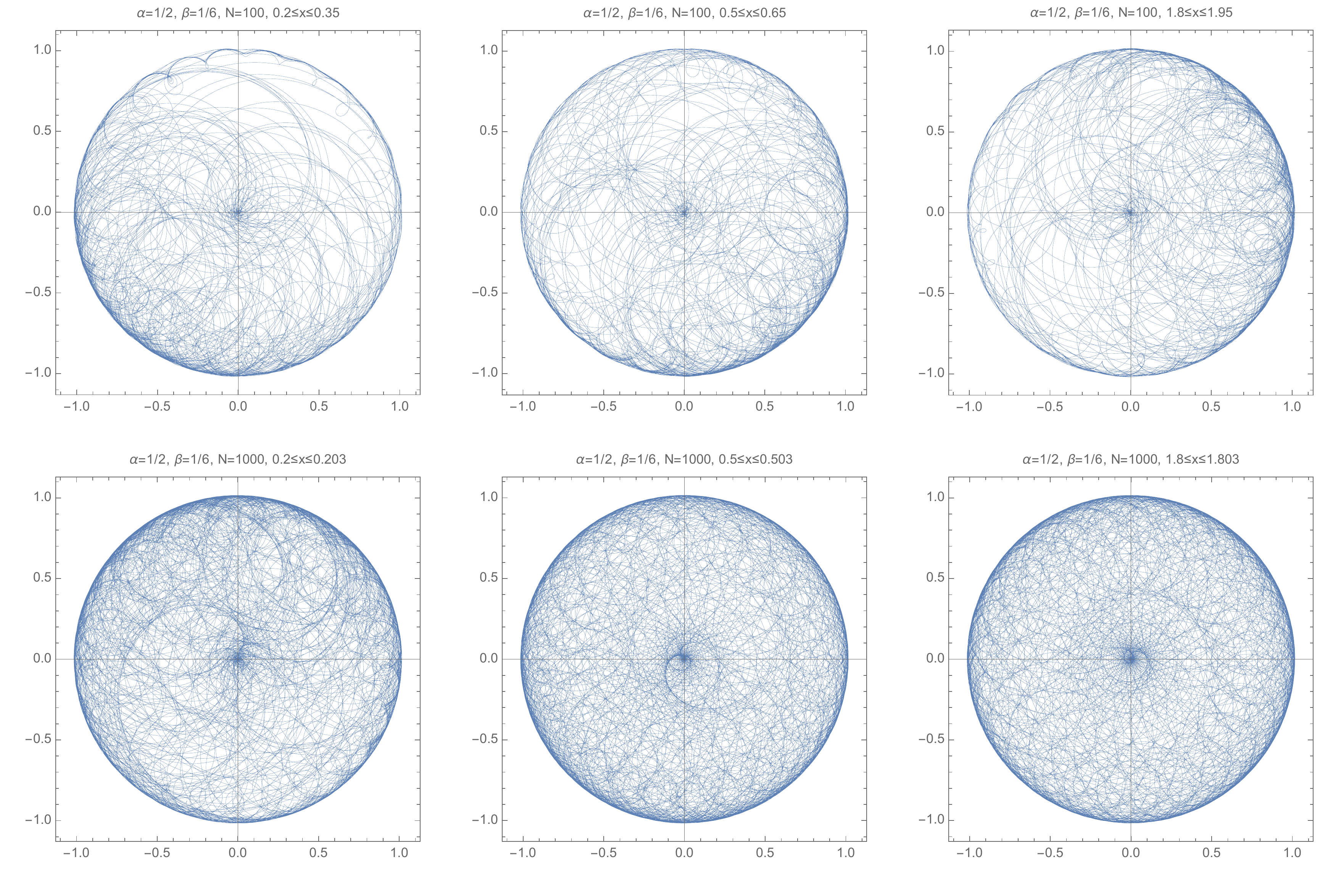}
\caption{Illustration of Theorem \ref{thm:jacobi-theta-bound-intro}. Curves $x\mapsto \varepsilon^{1/4}\vartheta_3(\frac{1}{2}+\frac{x}{6},x+i\varepsilon)\in\C$  with  $\varepsilon=N^{-2}$ and $N=100$ (top row) or $N=1000$ (bottom row). Each panel indicates the range of $x$.}
    \label{fig:six-figures-theta3-example}
\end{figure}

\begin{figure}[ht!]
    \centering
    \hspace{-.9cm}\includegraphics[width=16.6cm]{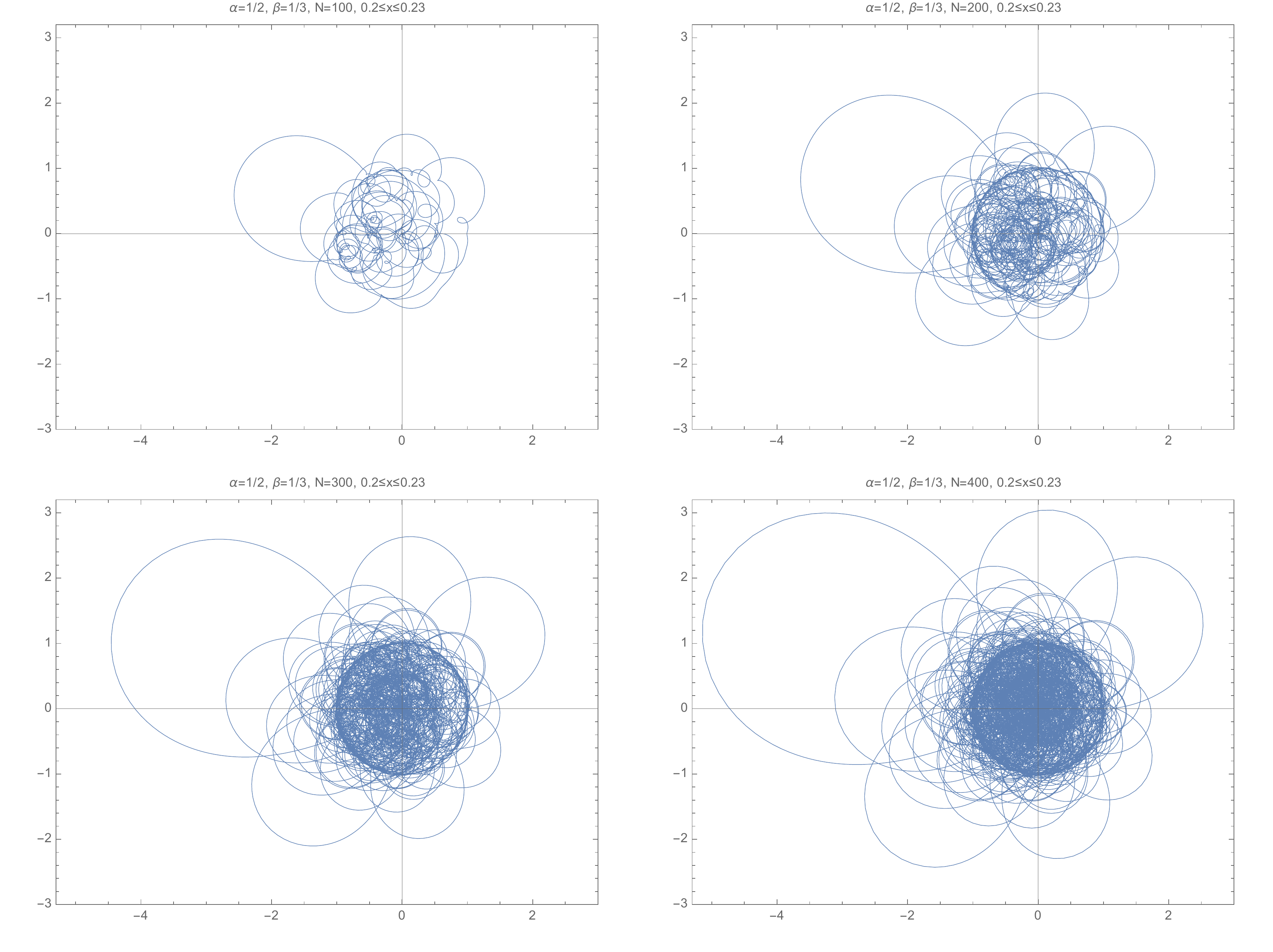}
    \caption{Curves $[0.2,0.23]\ni x\mapsto \varepsilon^{1/4}\vartheta_3(\frac{1}{2}+\frac{x}{3},x+i\varepsilon)\in\C$  with  $\varepsilon=N^{-2}$ and $N\in\{100,200,300,400\}$. The pair $(\frac{1}{2},\frac{1}{3})$ does not satisfy the hypotheses of Theorem \ref{thm:jacobi-theta-bound-intro}.}
    \label{fig:six-figures-theta3-nonexample}
\end{figure}

Figure \ref{fig:six-figures-theta3-example} illustrates Theorem \ref{thm:jacobi-theta-bound-intro} for $(\alpha,\beta)=(\frac{1}{2},\frac{1}{6})$ (that is $a=3,b=1,m=3$) for various values of $x$ and of $\varepsilon=N^{-2}$. For contrast, Figure \ref{fig:six-figures-theta3-nonexample} shows the lack of a uniform bound for $\varepsilon^{1/4}\vartheta_3(\frac{1}{2}+\frac{x}{3},x+i\varepsilon)$  as $\varepsilon=N^{-2}$ decreases. In this case the pair $(\alpha,\beta)=(\frac{1}{2},\frac{1}{3})$ does not satisfy the hypotheses of Theorem   \ref{thm:jacobi-theta-bound-intro} since $a=3, b=2, m=3$ are not all odd.

There are many results in the literature concerning upper bounds for generalised theta sums as defined in \eqref{def:S_N^f(x,alpha,beta)-intro}. These estimates have modern applications, for instance, in  \cite{Buterus-Gotze-Hille-Margulis}, \cite{Gotze-lattice-points}, \cite{Marklof2002Duke}, \cite{Marklof2003ann} they are used to understand the value distribution of quadratic forms.

 Estimates for $S_N := S_N^f$ where $f = \ind_{(0,1]}, k = 1$ have received much attention over the years. In the case where $x, \alpha, \beta \in \Q$, $S_N$ reduces to a quadratic Gauss sum for which various bounds are classical, e.g. if $\gcd (a,q) = 1$ and $N \leq q$ then $|S_N(\tfrac{a}{q}; 0, 0)| = O(\sqrt q)$. See \cite{Korolev-incomplete-gauss-sums}, \cite{Oskolov-gauss-sums} for further details.

The detailed study of $S_N$ for $x \in \R$ was initiated by Hardy and Littlewood in \cite{Hardy-Littlewood1914partI}, who were attracted by its ``interesting and beautiful properties''. In particular, they prove an approximate functional equation which they use to obtain various bounds for $S_N$, typically with some restriction on $x$, e.g. for $x$ of bounded type they prove that $S_N(x; \alpha,0) \ll_x \sqrt N$. 

For generic $x$ it is known that $|S_N(x; \alpha,\beta)| \ll_x \sqrt N \log N$, first in the case where $(\alpha, \beta) = (0,0)$ 
 by Fiedler, Jurkat and K\"{o}rner \cite{Fiedler-Jurkat-Korner77}, and then for every $(\alpha,\beta) \in \R^2$ by Flaminio and Forni \cite{Flaminio-Forni2006}, even with a reduction of the power of $\log N$, see Remark \ref{Remark-FF-bound}. Analogous bounds for theta sums in higher rank, including the  were obtained by Cosentino and Flaminio \cite{Cosentino-Flaminio}, with recent improvements by Marklof and Welsh \cite{Marklof-Welsh-higher-rank-theta-sums}, \cite{Marklof-Welsh-higher-rank-theta-sumsII}. For bounds on Weyl sums of arbitrary degree (not just theta sums, in which the degree of the polynomial in the exponential sum is $2$), see the recent work of Flaminio and Forni \cite{flaminio2023equidistribution} and the references therein.

In light of these results, we see that the behaviour of $|S^f_N(x;\vecalf,\vecbeta)|$ for $(\vecalf,\vecbeta)$ as in Theorem \ref{thm:main-thm-intro} and Corollary \ref{cor:main-cor-intro} is far from typical. In particular, the implied constant in \eqref{cor:main-cor-intro-bound} is independent of $x$. Therefore, for $(\alpha,\beta)$ as in the statement of Corollary \ref{cor:main-cor-intro} there exists $R_0 >0$ such that for any $R > R_0$
\begin{align}
     \Leb \{x \in \R : \tfrac{1}{\sqrt N} S^f_N (x; \alpha,\beta) > R\} = 0,
\end{align}
for every $N$. It follows that for $R > R_0$
\begin{align}
    \lim_{N\to \infty} \Leb \{x \in \R : \tfrac{1}{\sqrt N} S^f_N (x; \alpha,\beta) > R\} = 0,
\end{align}
and so the limiting distribution (as $N\to\infty$) of $\tfrac{1}{\sqrt N} S^f_N(x;\alpha,\beta)$ for random $x$ must be compactly supported in this case. In this way we see that Theorem \ref{thm:main-thm-intro} implies Theorem 1.0.7 (i) in \cite{Cellarosi-Osman-rational-tails}.

Our general approach aligns with that of \cite{Marklof-1999}, \cite{Marklof1999b}, \cite{Marklof2003ann} which interprets $\tfrac{1}{\sqrt N}|S^f_N|$ as an automorphic function $|\Theta_f|$ evaluated along special curves, known as horocycle lifts.

In Section \ref{background} we outline the construction of the (projective) Shr\"{o}dinger-Weil representation, a unitary representation of Jacobi group on $\LtRk$. We then use this representation in Section \ref{sec:Jacobi-theta-function} to define a special real-valued function $|\Theta_f|$ over the group $G := \sltr \ltimes \R^{2k}$, which is a subgroup of the Jacobi group, and $f: \R \to \R$ is some sufficiently regular weight function. 
The sum $\tfrac{1}{\sqrt N} |S^f_N|$ agrees with $|\Theta_f|$ along special curves $\calC\tovab_N$, known as horocycle lifts. Special invariance properties of $|\Theta_f|$ imply that it can be viewed as a function on a non-compact, finite volume homogeneous space $\GaG$.
In Section \ref{uniform-bound-section} we prove a uniform bound for $|\Theta_f|$, provided certain parameters avoid certain (explicit) regions of growth in $\GaG$. We then give a simple, explicit condition on $(\vecalf, \vecbeta)$ so that the curve  $\Gamma\calC\tovab_N$ (a horocycle lift viewed as a curve in $\GaG$) is bounded away from the regions of growth, uniformly in $N$. In Section \ref{study-of-orbits} we give explicit conditions on the pair $(\vecalf,\vecbeta)$ so that $\Gamma\calC\tovab_N$ avoids regions of growth, yielding Theorem \ref{thm:main-thm-intro}. 

The final Section \ref{section:applications} contains applications of our main theorem to produce a bound for the classical Jacobi theta function $\vartheta_3$ (Theorem \ref{thm:jacobi-theta-bound-intro}) and a bound of the form $|S_N(x;\alpha,\beta)|\ll_\varepsilon N^{1/2+\varepsilon}$ uniformly in $x$ and $N$ (Corollary \ref{cor:S_N-bound}) for the rational pairs $(\alpha,\beta)$ we consider. 


%% file: background_material-01.tex
\section{Representation-Theoretical Preliminaries}\label{background}
In this section we briefly introduce the necessary representation-theoretical ingredients needed to define (in Section \ref{sec:Jacobi-theta-function}) the Jacobi  theta function $\Theta_f$ for any sufficiently regular weight function $f$.  For further details, we refer the reader to  \cite{Marklof2003ann}.

\subsection{The Heisenberg Group its  Schr\"{o}dinger Representation}\label{schrodinger-repn-section}
Let $\omega:\R^{2k}\times\R^{2k}\to\R$ be the standard symplectic form,
\begin{align}\label{eq:symplectic-form-omega}
\omega\!\lp \ve{\vecxi_1^{\bianco}}{\vecxi_2^{\bianco}}, \ve{\vecxi'_1}{\vecxi'_2}\rp = \vecxi_1 \cdot \vecxi'_2 - \vecxi'_1\cdot\vecxi_2, 
\end{align}
where $\vecxi_1, \vecxi_2, \vecxi'_1, \vecxi'_2 \in \R^k.$ 
We define the Heisenberg group as $\Heik=\R^{2k} \times \R$ with multiplication law 
\begin{align}
  (\vecxi, \zeta)(\vecxi', \zeta') = (\vecxi + \vecxi', \zeta + \zeta' + \tha \omega(\vecxi, \vecxi')).
\end{align}
Note that $\R^{2k}$ is isomorphic to the quotient of $\Heik$ by its centre. 
The \emph{Schr\"{o}dinger representation} $W$ is a representation of $\Heik$ on  $\mathcal{U}(\LtRk)$, the group of unitary operators on $\LtRk$. Specifically, using the fact that each element in $\Heik$ can be decomposed as
\begin{align}\label{Heis-decomp}
  \veH{\vecxi_1}{\vecxi_2}{\zeta} = \veH{\vecxi_1}{\bn}{0} \veH{\bn}{\vecxi_2}{0} \veH{\bn}{\bn}{\zeta - \frac{\vecxi_1 \cdot\vecxi_2}{2}},
\end{align}
we define $W$ via  
\begin{align}
  &\ls W\!\sveH{\vecxi_1}{\bn}{0} f\rs (\vecw) = \e{(\vecxi_1 \cdot \vecw}f(\vecw)\label{schrodinger-100},\\
  &\ls W\!\sveH{\bn}{\vecxi_2}{0} f\rs (\vecw) = f(\vecw - \vecxi_2)\label{schrodinger-010},\\
  &\ls W\!\sveH{\bn}{\bn}{\zeta} f\rs (\vecw) = \e{\zeta} f(\vecw)\label{schrodinger-001},
\end{align}
where $f\in\LtRk$. For further details on this construction, see Section 1.2 of \cite{Lion-Vergne}.

\subsection{$\mathrm{Sp}(k,\R)$ and its Projective Shale-Weil Representation}\label{shale-weil-repn-section}
For any $M \in \mathrm{Sp}(k,\R)$ we may define a new representation of $\Heik$ as $W_M (\vecxi, \zeta) := W(M\vecxi, \zeta)$ where $W$ is the Schr\"{o}dinger representation defined in \eqref{schrodinger-100}-\eqref{schrodinger-001}. By the Stone-von Neumann theorem 
any such representation is irreducible, and unitarily equivalent. Therefore, there exists a unitary operator $R(M)$ on $\LtRk$ such that
\begin{align}
    R(M)W(\vecxi, \zeta)R(M)^{-1} = W_M(\vecxi,\zeta). \label{Shale-Weil-defn}
\end{align}
By Schur's Lemma the map $M \mapsto R(M)$ is unique up to a phase, that is
\begin{align}
  R(M_1M_2) = c(M_1, M_2)R(M_1)R(M_2) \label{cocycle-relation}
\end{align}
where $|c(M_1,M_2)| = 1$. The map $R: \mathrm{Sp}(k,\R) \to \mathcal{U}(\LtRk)$ is therefore a projective unitary representation of $\sltr$, known as the \emph{projective Shale-Weil representation}. The phase $c(M_1,M_2)$ can be computed explicitly (see Theorem 1.6.11 in \cite{Lion-Vergne}), but we shall not need it in what follows.

\subsection{The Jacobi Group and its Projective Schr\"odinger-Weil Representation}
We define the Jacobi group as $\mathrm{Sp}(k,\R) \ltimes \Heik$ whose group law is
\begin{align}
  (M;\vecxi,\zeta)(M';\vecxi', \zeta') = (MM' ; \vecxi + M \vecxi', \zeta + \zeta' - \tha \omega (\vecxi, M \vecxi')),
\end{align}
where $\omega$ is as in \eqref{eq:symplectic-form-omega}. 
Using the Schr\"odinger representation $W$ of $\Heik$ from Section \ref{schrodinger-repn-section} and the projective Shale-Weil representation $R$ of $\sltr$ from Section \ref{shale-weil-repn-section}, we can define the \emph{projective Schr\"odinger-Weil representation} of the Jacobi group as
\begin{align}
  R(M; \vecxi,\zeta) := W(\vecxi,\zeta)R(M). \label{schrodinger-weil-repn-defn}
\end{align}

\section{The Jacobi Theta Function}\label{sec:Jacobi-theta-function}

Observe that we may embed $\sltr$ in $\mathrm{Sp}(k, \R)$ via $\smatr{a}{b}{c}{d} \mapsto \smatr{aI_k}{bI_k}{cI_k}{dI_k}$ where $I_k$ is the $k\times k$ identity matrix. Given a  function $f: \R^k \to \R$ we  define (up to a phase) a theta function $\Theta_f:\sltr \ltimes \Heik \to\C$ given by
\begin{align}
  \Theta_f(M; \vecxi, \zeta) := \sum_{\vecn \in \Z^k} [R(M; \vecxi, \zeta)f](\vecn), 
  \label{Jacobi-theta-sum-1}
\end{align}
provided the series in \eqref{Jacobi-theta-sum-1} converges absolutely (we shall make sufficient assumptions for this to hold, see Section \ref{sec:regular-weights}). The function is defined up to a phase because $R(M)$ in \eqref{Shale-Weil-defn}, and therefore $R(M; \vecxi, \zeta)$ in \eqref{schrodinger-weil-repn-defn}, are defined up to a phase. To properly define $\Theta_f$, one can pass to the universal cover $\tsltr\ltimes\Heik$ as done in \cite{Marklof-1999}, \cite{Marklof2007b}, \cite{Cellarosi-Marklof}, \cite{Cellarosi-Osman-rational-tails} when $k=1$. However, since our aim is to study $|\Theta_f|$, this definition of $\Theta_f$ will suffice. Let us see that  $|\Theta_f|$ may be viewed as a real-valued function on the group
  $G:= \sltr \ltimes \R^{2k}$,
with group law
\begin{align}
  (M; \vecxi)(M';\vecxi') = (MM' ; \vecxi + M \vecxi'). \label{spkr-r2k-law}
\end{align}
In fact, by $\eqref{schrodinger-weil-repn-defn}$, we see that 
  $R(M; \vecxi, \zeta) = W(\bn, \zeta)W(\vecxi, 0)R(M)$.
Therefore, by  \eqref{schrodinger-001} and \eqref{Jacobi-theta-sum-1}, it follows that
\begin{align}
  |\Theta_f(M; \vecxi, \zeta)| = \labs e(\zeta) \sum_{n \in \Z^{k}} R(M; \vecxi, 0) f(\vecn)\rabs = |\Theta_f(M; \vecxi, 0)|,
\end{align}
and so we may define $|\Theta_f|:G\to\R$ as
\begin{align}
  |\Theta_f| (M;\vecxi) := |\Theta_f (M;\vecxi, 0)|.
\end{align}

\subsection{Writing $|\Theta_f|$ in coordinates on $\sltr \ltimes \R^{2k}$}
Through the embedding of $\sltr$ in $\mathrm{Sp(k, \R)}$, we obtain an action of $\sltr$ on $\R^{2k}$ coming from the group law \eqref{spkr-r2k-law} of $\mathrm{Sp}(k,\R) \ltimes \R^{2k}$, namely
\begin{align}
  \matr{a}{b}{c}{d} \cdot \ve{\vecxi_1}{\vecxi_2} = \ve{a \vecxi_1 + b\vecxi_2}{c\vecxi_1 + d\vecxi_2}.\label{SL2-action-on-R2k}
\end{align}
By Iwasawa decomposition, any $M \in \sltr$ may be written as
\begin{align}
  M = n_x a_y k_{\phi} \label{Iwasawa-decomp}
\end{align}
where $x + iy \in \h := \{z \in \C : \Im (z) > 0\}$, $\phi \in [0, 2\pi)$, and 
\begin{align}
  n_x := \matr{1}{x}{0}{1},\:\:\:\: a_y := \matr{y^{1/2}}{0}{0}{y^{-1/2}}, \:\:\:\: k_{\phi} := \matr{\cos \phi}{-\sin \phi}{\sin \phi}{\cos \phi}.
\end{align}
It can be shown (see, e.g., Section 1.6 of \cite{Lion-Vergne}) that
\begin{align}
  [R(n_x)f](\vecw) &= \e{\tha \|\vecw\|^2 x}f(\vecw),\label{horocycle-action-on-L2}\\
  [R(a_y)f](\vecw) &= y^{\frac{k}{4}}f\!\lp y^{\frac{1}{2}}\vecw\rp,\label{geodesic-action-on-L2}\\
  [R(k_{\phi})f](\vecw) &= 
  \begin{dcases*} 
    f(\vecw) & if $\phi = 0 \!\!\!\pmod {2\pi}$,\label{rotation-action-on-L2}\\
    f(-\vecw) & if $\phi = \pi \!\!\!\pmod {2\pi}$,\\
  |\sin \phi|^{-\frac{k}{2}}\int_{\R^k} \e{\tfrac{\ha(\|\vecw\|^2 + \|\vecv\|^2)\cos \phi - \vecw\cdot\vecv}{\sin \phi}}f(\vecv) \:\de \vecv & if $\phi \neq 0\!\!\!\pmod{\pi}$.
  \end{dcases*}
\end{align}
In view of \eqref{Iwasawa-decomp}, formul\ae\: \eqref{horocycle-action-on-L2}--\eqref{rotation-action-on-L2} are enough to define the projective Shale-Weil representation restricted to $\sltr$. 
Define 
\begin{align}
  f_\phi := R(k_{\phi})f. \label{phi-transform-definition}
\end{align}
Identifying $(M,\vecxi)\in\sltr\ltimes\R^{2k}$ with $(z,\phi ;\vecxi) \in \h \times [0,2\pi) \times \R^{2k}$ via \eqref{Iwasawa-decomp},  the modulus of the theta function  $|\Theta_f|$ may be written in `Iwasawa coordinates' as
\begin{align}\label{Jacobi-theta-sum-2}
  |\Theta_f (z, \phi ; \vecxi)| = y^{k/4}\labs\sum_{\vecn \in \Z^{k}}f_{\phi}((\vecn - \vecxi_2)y^{1/2})\e{\tha\|\vecn - \vecxi_2\|^2 x + n \vecxi_1}\rabs,
\end{align}
We now see how the generalised quadratic Weyl sums \eqref{def:S_N^f(x,alpha,beta)-intro} can be related to certain values of $|\Theta_f|$. By choosing $(z, \phi; \vecxi) = ( x + iN^{-2}, 0; \sve{\vecalf + \vecbeta x}{0})$ we have
\begin{align}
  |\Theta_f (x + \tfrac{i}{N^2}, 0 ; \sve{\vecalf + \vecbeta x}{0})| &= N^{-k/2}\!\sum_{\vecn \in \Z^k} f(\tfrac{1}{N}\vecn)\,\e{(\tha \|\vecn\|^2 + \vecbeta \cdot \vecn)x + \vecalf \cdot\vecn} \nonumber\\&= N^{-k/2}S^f_N (x; \vecalf, \vecbeta).\label{key-relation-theta-sum-to-function}
\end{align}
Hence, in light of \eqref{key-relation-theta-sum-to-function}, bounds of the form $|S_N^f(x,\vecalf,\vecbeta)|\ll N^{k/2}$ can be obtained by bounds $|\Theta_f (x + \tfrac{i}{N^2}, 0 ; \sve{\vecalf + \vecbeta x}{0})|\ll 1$.

\subsection{Regular Cut-off Functions}\label{sec:regular-weights}
When we introduced $S_N^f$ in \eqref{def:S_N^f(x,alpha,beta)-intro}, we assumed that $f$ is a Schwartz function. We now relax this assumption. It is apparent from \eqref{Jacobi-theta-sum-2} that $|\Theta_f|:\sltr \ltimes \R^{2k} \to \R$ is not necessarily well defined pointwise for arbitrary $f$. 
For instance, as pointed out in Section 2.6 of \cite{Cellarosi-Marklof}, taking $k = 1$ and $f = \ind_{(0,1]}$ we see that  $\Theta_f(z, 0; \vecxi, \zeta)$ converges  since it is a finite sum, but $f_{\pi/2} (w)$ decays too slowly as $|w| \to \infty$, and hence the series defining $\Theta_f(z, \pi/2; \vecxi, \zeta)$ does not converge absolutely.
Define
\begin{align}
  \kappa_\eta (f) := \sup_{\vecw, \phi} (1 + \|\vecw\|^2)^{\eta/2}|f_\phi (\vecw)|,\label{kappa-norm-def}
\end{align}
and consider the function space
\begin{align}
  \Si_\eta (\R^k) := \{f: \R^k \to \R : \kappa_\eta (f) < \infty\}.\label{def-S_eta}
\end{align}
We say that a cut-off function $f: \R^{k} \to \R$ is \emph{regular} if it belongs to $\Si_\eta(\R^k)$ for some $\eta>k$. Note that the regularity of $f$ is sufficient to 
guarantee that the series \eqref{Jacobi-theta-sum-2} defining $|\Theta_f|$ is absolutely convergent for every $(z, \phi; \vecxi) \in \h \times [0,2\pi) \times \R^{2k}$.
Regular cut-off functions generalise Schwartz functions. In fact, it can be shown that Schwartz functions belong to $\Si_\eta (\R^k)$ for every $\eta>1$, 
  see Lemma 4.3 in \cite{Marklof2003ann}. 
The following lemma 
will be useful in the proof of the main theorem.

\begin{lemma}\label{kappa_shift}
  Let $h := \sveH{\veca}{\vecb}{c}\in \Heik$. If $f\in \Si_{\eta}(\R^k)$ for some $\eta > 1$ then $W(h)f\in\Si_\eta(\R^k)$.
  More precisely, 
  \begin{align}
      \kappa_{\eta}(W(h)f) \leq 3^\eta (1 + \|\veca\|^2 + \|\vecb\|^2 )^{\eta/2}\, \kappa_{\eta}(f).
  \end{align}
\end{lemma}

\begin{proof}
  Due to \eqref{Heis-decomp} and \eqref{schrodinger-001}, it is enough to consider $h=(\sve{\veca}{\bn}, 0)(\sve{\bn}{\vecb},0)$. 
 Using \eqref{Shale-Weil-defn} and the fact that  $W_{k_\phi}(h)=\ W\!\lp \sve{\veca\cos \phi - \vecb\sin \phi}{\veca\sin \phi+ \vecb \cos \phi}, \tfrac{\veca\cdot\vecb}{2}\rp$,  we get
  \begin{align}
      |[W(h) f]_{\phi}(\vecw)| = \labs \left[R(k_\phi)W\!\lp\sve{\veca}{\vecb}, 0\rp  f\right]\!(\vecw)\rabs &= \labs W\!\lp \sve{\veca\cos \phi - \vecb\sin \phi}{\veca\sin \phi + \vecb \cos \phi}, 0\rp \cdot f_{\phi} (\vecw)\rabs \\ &= |f_{\phi}(\vecw - (\veca\cos \phi + \vecb \sin \phi) )|.
  \end{align}

  Now, note that
  \begin{align}
      \kappa_{\eta}(h\cdot f) &= \sup_{\vecw,\phi} \:(1 + \|\vecw\|^2)^{\eta/2}|f_{\phi}(\vecw - (\veca\cos \phi + \vecb \sin \phi) )|\\
      &= \sup_{\vecw,\phi}\: (1 + \|\vecw + (\veca\cos \phi + \vecb \sin \phi)\|^2)^{\eta/2}|f_{\phi}(\vecw)|. \label{change-of-variables-in-sup}
  \end{align}
  Applying the triangle inequality and the AM-GM inequality to the multiplier in \eqref{change-of-variables-in-sup} we have
  \begin{align}
    &\sup_{\vecw,\phi}\: (1 + \|\vecw + (\veca\cos \phi + \vecb \sin \phi)\|^2)^{\eta/2}|f_{\phi}(\vecw)|\\
    &\leq 3^{\eta} \sup_{\vecw, \phi}\: ((1 + \|\veca\|^2 + \|\vecb\|^2) +  \|\vecw\|^2 )^{\eta/2}|f_{\phi}(\vecw)|\\
    &\leq 3^\eta(1+ \|\veca\|^2 + \|\vecb\|^2)^{\eta/2}\kappa_\eta (f),
  \end{align}
  which proves the lemma.
\end{proof}

\subsection{Invariance Properties}\label{section-invariance-properties}
Using \eqref{Iwasawa-decomp} to parameterise $\sltr$ by $\h \times [0, 2\pi)$, we obtain the following  transitive action of $\sltr$ on $\h \times [0, 2\pi)$:
\begin{align}\label{unit-tangent-bundle-action}
\smatr{a}{b}{c}{d}. (z, \phi) = \lp \frac{az + b}{cz + d}, \phi + \arg (cz + d)\rp.
\end{align}
We then define an action of special affine linear group $G$ with group law 
on $\h \times [0,2\pi) \times \R^{2k}$ as 
\begin{align}\label{action-of-G-on-h-times-R2}
\lp \smatr{a}{b}{c}{d};\vecv\rp\!. (z, \phi;\vecxi) = \lp \frac{az + b}{cz + d}, \phi + \arg (cz + d); \vecv+\smatr{a}{b}{c}{d} \cdot \vecxi \rp.
\end{align}
where $\smatr{a}{b}{c}{d} \cdot \vecxi$ is as in \eqref{SL2-action-on-R2k}
The following transformation formul\ae\: for $\Theta_f$ follow from Poisson summation and the definition of $\Theta_f$. 
\begin{theorem}[\cite{Marklof2003ann}, ``Jacobi 1-3"]\label{theorem-Jacobi-1-2-3}
  Let $f \in \Si_\eta(\R^k)$, with $\eta > k$. Let $\Delta<G$ be the group generated by
  \begin{align}
    \rho_1 := \lp\smatr{0}{-1}{1}{0}; \bn\rp,\:\:\rho_2 := \lp \smatr{1}{1}{0}{1}; \sve{\vecs}{\bn}\rp,\:\: \rho_{\vecm, \vecn} := \lp I; \sve{\vecm}{\vecn}\rp,
  \end{align}
  where $\vecs = (\tha,\dots, \tha)^\top$, and $\vecm, \vecn \in \Z^{k}$. 
  Then
  \begin{align}
   |\Theta_f(\rho \cdot (z,\phi; \vecxi))| = |\Theta_f(z, \phi; \vecxi)|
  \end{align}
  for any $\rho \in \Delta$.
\end{theorem}

It follows that for $f\in \Si_{\eta}(\R^k)$ with $\eta > 1$, 
the function $|\Theta_f|$ is a well-defined real-valued function on the group $\DaG$. 
For convenience, we pass to a finite index subgroup $\Gamma$ of $\Delta$, whose generators are 
\begin{align}
\gamma_1 = \rho_1, \: \gamma_2 = \rho^2_1\rho_{\vecu, \bn}=\lp\smatr{1}{2}{0}{1}; \bo{0}\rp, \: \gamma_{\vecm,\vecn} = \rho_{\vecm, \vecn}, 
\end{align}
 where $\vecu= (-1,\dots, -1)^\top$, and $\vecm, \vecn \in \Z^k$.

Clearly, $\Gamma = \Gamma_\theta \ltimes \Z^{2k}$, where the lattice 
\begin{align}
\Gamma_\theta = \langle \smatr{0}{-1}{1}{0}, \smatr{1}{2}{0}{1}\rangle<\sltz\label{def-Gamma_theta}
\end{align} is the so-called \emph{theta group}.
%
We therefore view $|\Theta_f|$ as a function on $\GaG$ or $\Gamma \backslash (\h \times [0,2\pi) \times \R^{2k})$, whichever is more convenient. This quotient can be shown to be non-compact and of with finite volume according to the Haar measure on $G$.

\subsection{Fundamental Domains}\label{section-fundamental-domains}

As the action of $G$ on the space $\h \times [0,2\pi) \times \R^{2k}$ defined in \eqref{action-of-G-on-h-times-R2} is transitive, we may identify $\GaG$ and $\Gamma \backslash (\h \times [0,2\pi) \times \R^{2k})$.  
\begin{proposition}
  A fundamental domain for the action of $\Gamma$ on $\h \times [0,2\pi) \times \R^{2k}$ is given by
  \begin{align}\label{fund-dom-gamma}
    \F_{\Gamma} = \F_{\Gamma_\theta} \times [0,\pi) \times [0,1)^{2k},
  \end{align}
   where $\F_{\Gamma_\theta} = \{z \in \h : |z| > 1, |z - 2| > 1, 0 \leq \Re z < 2\}$ is a fundamental domain for the action of $\Gamma_\theta$ on $\h$.
\end{proposition}
\begin{proof}
The proof is essentially identical to the proof of Proposition 2.5.1 in \cite{Cellarosi-Osman-rational-tails}. 
\end{proof}

The choice of  fundamental domain $\F_{\Gamma_\theta}$ has two cusps, one at $i\infty$ of width two and another at $1$ of width one. 
For reasons that will become apparent later, we split the fundamental domain $\F_{\Gamma}$ into $\Finfty$ and $\Fone$ in order to isolate these cusps, namely
\begin{align}
  &\Finfty := \{(z, \phi; \vecxi, \zeta) \in \F_{\Gamma} : |z| \geq 1\} \label{F-infty-def},\\
  &\Fone := \{(z, \phi; \vecxi, \zeta) \in \F_{\Gamma} : |z| < 1\}\label{F-one-def}.
\end{align}

%% file: uniform-bounds-v1.tex
\section{Uniform Bounds}\label{uniform-bound-section}
\subsection{Uniform Bounds for fixed $\vecxi$}
We begin with a lemma along the lines of Lemma 2.1 in \cite{Cellarosi-Marklof}.
Define, for $\vecv\in\R^k$,
\begin{align}
    \theta_k(\vecv) := \min_{\vecn\in \Z^k} \|\vecn - \vecv\|, \label{closest-integer-vector}
\end{align}
i.e. Euclidean distance to the closest integer point.
\begin{lemma}\label{L2.1Bdd}
    Suppose $\vecxi = \sve{\vecxi_1}{\vecxi_2} \in \R^{2k}$ is such that $\theta_k(\vecxi_2) > 0$. Let $f \in \Si_{\eta}(\R^k)$ with  $\eta > k$. Then there exits a constant $C = C(k,\vecxi_2,\eta)>0$ such that 
    \begin{align}
        |\Theta_f (x + iy, \phi, \vecxi)| \leq C\, \kappa_\eta (f).
    \end{align}
    for every $y \geq \tha$, and every $x, \phi$ and $\vecxi_1$.
\end{lemma}

\begin{proof}
    As  $\vecxi_2 \not\in \Z^k$ it follows that 
    \begin{align}
        |\Theta_f (x + iy, \phi, \vecxi)| &\leq y^{k/4}\sum_{\vecn \in \Z^k} |f_{\phi}((\vecn-\vecxi_2)y^{1/2})|\\
        &\leq \kappa_{\eta}(f)\, y^{k/4} \sum_{\vecn \in \Z^k}\frac{1}{(1 + \|\vecn - \vecxi_2\|^2y)^{\eta/2}}\\
        &\leq \kappa_{\eta}(f)\, y^{-\frac{2\eta - k}{4}} \sum_{\vecn \in \Z^k} \frac{1}{\|\vecn - \vecxi_2\|^{\eta}}.
    \end{align}
    As $y \geq \tha$ it follows that $y^{-\frac{2\eta-k}{4}} \leq 2^{\frac{2\eta - k}{4}}$ since $\eta > k$. Let 
    \begin{align}
    C(k,\vecxi_2,\eta) := 2^{\frac{2\eta - k}{4}}\sum_{\vecn \in \Z^k} \frac{1}{\|\vecn - \vecxi_2\|^{\eta}},
    \label{def-C(k,xi2,eta)}
    \end{align}
   which is a convergent series for $\eta > k$, as $\theta_k(\vecxi_2) > 0$.

\end{proof}
\begin{remark}\label{remark-size-constant}
    Let $\mathcal M := \{ \vecn \in \Z^k : |\vecn - \vecxi_2| = \theta_k (\vecxi_2)\}$. Note that $1\leq |\mathcal M| \leq 2^k$. The series 
   \begin{align}
     \sum_{\vecn \in \Z^k} \frac{1}{\|\vecn - \vecxi_2\|^{\eta}} = \frac{|\mathcal M|}{\theta_k(\vecxi_2)^\eta} + O_\eta (1),
    \end{align}
    and so the closer $\vecxi_2$ is to an integer, the larger the constant $C(k,\vecxi_2,\eta)$ in \eqref{def-C(k,xi2,eta)} becomes. 
\end{remark}

\begin{cor}\label{corFinfty}
    Suppose that $(x + iy; \vecxi) \in \Finfty$ where $\vecxi = \sve{\vecxi_1}{\vecxi_2} \in \R^{2k}$ is such that $\theta_k(\vecxi_2) > 0$. Let $f\in \Si_\eta(\R^k)$ with $\eta > k$. Then 
    \begin{align}
        |\Theta_f (x + iy, \phi, \vecxi)| \leq C(k, \vecxi_2, \eta)\, \kappa_{\eta}(f),
    \end{align}
    where $C$ is as in Lemma \ref{L2.1Bdd}.
\end{cor}
 
\begin{proof}
    This is immediate from Lemma \ref{L2.1Bdd} 
    since $y \geq \tha$ for any $(x + iy; \vecxi) \in \Finfty$. 
\end{proof}
\begin{cor}\label{corFone}
    Suppose $(x + iy; \vecxi) \in \Fone$ and $\vecxi = \sve{\vecxi_1}{\vecxi_2} \in \R^{2k}$. Let  $\vecs$ be as in Theorem \ref{theorem-Jacobi-1-2-3}, define $\hat\vecxi_2 := \vecxi_2 - \vecxi_1 + \vecs$, and suppose $\theta_k(\hat\vecxi_2) > 0$. Let $f\in \Si_\eta(\R^k)$ with $\eta > k$. Then 
    \begin{align}
        |\Theta_f (x + iy, \phi, \vecxi)| \leq C(k,\hat\vecxi_2, \eta)\,\kappa_{\eta}(f),\label{statement_corFone}
    \end{align}
    where $C$ is as in Lemma \ref{L2.1Bdd}.
\end{cor}

\begin{proof}
    Let $\lp x + iy; \sve{\vecxi_1}{\vecxi_2}\rp \in \Fone$. Define $\rho := \lp\smatr{0}{1}{-1}{1};\sve{\bn}{\vecs}\rp = \rho_2^{-1}\rho_1^{-1}$ and set 
    \begin{align}
        \lp x' + iy'; \sve{\vecxi_1'}{\vecxi'_2}\rp = \rho\cdot \lp x + iy; \sve{\vecxi_1}{\vecxi_2}\rp.\label{g'_in_terms_of_rho.g}
    \end{align}
    A computation gives that $y' \geq \ha$ and $\vecxi'_2 =\hat\vecxi_2$. 
    Since $\hat\vecxi_2 \not\in \Z^{k}$,  it follows from Lemma \ref{L2.1Bdd} that
    \begin{align}
        \labs\Theta_f\!\lp x' + iy'; \sve{\vecxi_1'}{\vecxi'_2}\rp\rabs \leq C(k, \hat\vecxi_2,\eta)\,\kappa_{\eta}(f).\label{estimate_corFone_est'}
    \end{align}
    As $|\Theta_f|$ is invariant under $\rho \in \Delta$ (see Theorem \ref{theorem-Jacobi-1-2-3}), using \eqref{g'_in_terms_of_rho.g} we obtain \eqref{statement_corFone} from \eqref{estimate_corFone_est'}.
\end{proof}




%% file: key-ingredients.tex
\subsection{Bounds over  $\Gamma$-orbits 
}\label{crucial-relationship-section}

Recall that $\Gamma = \Gamma_{\theta} \ltimes \Z^{2k}$ acts on $\R^{2k}$ via the group law as follows:
\begin{align}
    \lp\smatr{a}{b}{c}{d} ; \sve{\vecm}{\vecn}\rp \cdot \sve{\vecxi_1}{\vecxi_2} = \sve{a \vecxi_1 + b \vecxi_2}{c \vecxi_1 + d\vecxi_2} + \sve{\vecm}{\vecn}.
\end{align}
This descends to an action on the torus $\T^{2k}$, or equivalently on $[0,1)^{2k} \bmod \Z^{2k}$. Set $\Orb_\Gamma \sve{\vecxi_1}{\vecxi_2}$ to be the orbit of $\sve{\vecxi_1}{\vecxi_2} \bmod \Z^{2k}$ under the action of $\Gamma$. For example, when $k=1$ we have $\Orb_{\Gamma}(\frac{1}{6},\frac{1}{6})=\{(\frac{1}{6}, \frac{1}{6}), (\frac{1}{6},\frac{1}{2}), (\frac{1}{6}, \frac{5}{6}), (\frac{1}{2}, \frac{1}{6}), (\frac{1}{2}, \frac{5}{6}), (\frac{5}{6}, \frac{1}{6}), (\frac{5}{6}, \frac{1}{2}), (\frac{5}{6}, \frac{5}{6})\}$, see also Figure \ref{fig:many-theta3-same-bound}

The following are one-parameter subgroups of $G$, which shall be referred to as the \emph{geodesic} and \emph{horocycle} flows, respectively:
\begin{align}
    \Phi &:= \lcur \Phi^t := \lp\smatr{e^{-t/2}}{0}{0}{e^{t/2}}; \bo{0}\rp \!:\: t \in \R\rcur,\label{geodesic-def}\\
    \Psi &:= \lcur \Psi^x := \lp \smatr{1}{u}{0}{1}; \bn\rp\::\: x \in \R\rcur\label{horocycle-def}.
\end{align}
We may rewrite the key relationship between $S^f_N$ and $\Theta_f$ given in \eqref{key-relation-theta-sum-to-function} using the right action of the geodesic and horocycle flows on the coset space $\GaG$ as 
\begin{align}
  N^{-k/2}|S^f_N (x; \vecalf,\vecbeta)|= |\Theta_f (x + iN^{-2}; \sve{\vecalf + \vecbeta x}{\bn})| = \labs\Theta_f \!\lp \Gamma (I; \sve{\vecalf + \vecbeta x}{\bn}) \Psi^x\Phi^{2\log N}\rp\rabs.\label{SfN-to-Theta-f}
\end{align}  
This observation first appeared in the works of Marklof \cite{Marklof-1999}, \cite{Marklof1999b} and plays a crucial role in the study of the limiting distribution of theta sums (i.e. when the weight function $f = \ind_{(0,1]}$). 

Recall the definition of $\theta_k$ as in \eqref{closest-integer-vector} and $\vecs$ as in Theorem \ref{theorem-Jacobi-1-2-3}. Define 
\begin{align}
    m_k(\vecalf,\vecbeta) := \inf \lcur\theta_k(\vecxi_2)\land \theta_k(\vecxi_2 - \vecxi_1 + \vecs) : \sve{\vecxi_1}{\vecxi_2} \in \Orb_\Gamma\sve{\vecalf}{\vecbeta}\rcur . \label{m-function}
\end{align}

\begin{theorem}\label{bounded-sums-chapter-5}
    Let $f \in \Si_\eta(\R^k)$ with $\eta > k$. Suppose $\sve{\vecalf}{\vecbeta} \in \R^{2k}$ with $m_k = m_k(\vecalf,-\vecbeta) > 0$. Then
    \begin{align}\label{bounded-sums-chapter-5-statement}
        |S_N^f(x; \vecalf,\vecbeta)| \ll_{k,m_k,\eta,\vecbeta,f} N^{k/2}
    \end{align}
    for every $x \in \R$ and every $N \in \N$.
\end{theorem}

\begin{proof}
    We have that
    \begin{align}
        \Gamma \lp I; \sve{\vecalf + \vecbeta x}{0}\rp \Psi^x\Phi^t = \Gamma \lp I; \sve{\vecalf}{-\vecbeta}\rp \Psi^x (I; \sve{0}{\vecbeta})\Phi^t = \Gamma \lp I; \sve{\vecalf}{-\vecbeta}\rp \Psi^x\Phi^t (I; \sve{\bn}{e^{-t/2}\vecbeta})
    \end{align}
    and so, by \eqref{Jacobi-theta-sum-1} 
    \begin{align}
        \labs\Theta_f \!\lp \Gamma (I; \sve{\vecalf + \vecbeta x}{0}) \Psi^x\Phi^{2\log N}\rp\rabs = |\Theta_{f_N} (\Gamma(I; \sve{\vecalf}{-\vecbeta} \Psi^x \Phi^{2\log N}))|, \label{f-to-f_N}
    \end{align}
    where $f_N = W(\sve{\bn}{N^{-1}\vecbeta},0)f$. 
    By \eqref{SfN-to-Theta-f} and \eqref{f-to-f_N} we have
    \begin{align}\label{4-19}
        N^{-k/2}|S_N^f (x; \vecalf,\vecbeta)| 
        = \labs\Theta_{f_N} \!(x + iN^{-2}, 0 ; \sve{\vecalf}{-\vecbeta})\rabs.
    \end{align}
    Let $\gamma \in \Gamma$ be the unique element such that $\gamma \cdot (x + iN^{-2}, 0 ; \sve{\vecalf}{-\vecbeta}) = \lp x' + iy' ,\phi'; \sve{\vecalf'}{\vecbeta'}\rp \in \calF_\Gamma$. Then, by the $\Gamma$-invariance of $|\Theta_{f_N}|$, we have
    \begin{align}\label{4-20}
        N^{-k/2}|S_N^f(x;\alpha,\beta)| = \labs\Theta_{f_N} \!\lp x' + iy' ,\phi'; \sve{\vecalf'}{\vecbeta'}\rp\rabs.
    \end{align}
    Let us apply either Corollary \ref{corFinfty} or Corollary \ref{corFone} depending on whether $\lp x' + iy' ,\phi'; \sve{\vecalf'}{\vecbeta'}\rp$ belongs to $\Finfty$ or $\Fone$.
    In the first case we have
    \begin{align}
        |\Theta_f (x' + iy' ,\phi'; \sve{\vecalf'}{\vecbeta'})| \leq C(k,\vecbeta',\eta)\, \kappa_{\eta}(f_N),
    \end{align}
    while in the second case we have the bound
    \begin{align}
        |\Theta_f (x' + iy' ,\phi'; \sve{\vecalf'}{\vecbeta'})| \leq C(k,\hat{\vecbeta'},\eta)\,\kappa_{\eta}(f_N),
    \end{align}
    where $\hat{\vecbeta'}=\vecbeta'-\vecalf'+\vecs$.
The previous two estimates imply
\begin{align}\label{bound-with-uniform-constant-over-orbit}
    |\Theta_f (x' + iy' ,\phi'; \sve{\vecalf'}{\vecbeta'})|\leq M\,\kappa_\eta(f_N),
\end{align}
where 
\begin{align}
    \label{constant-M}M=M(k,m_k,\eta)=\sup\lcur\max\{C(k,\vecbeta',\eta),C(k,\hat{\vecbeta'},\eta)\}:\:\sve{\vecalf'}{\vecbeta'}\in\Orb_\Gamma\sve{\vecalf}{-\vecbeta}\rcur.
\end{align}
    The constant $M$ is finite because of the assumption the $m_k(\vecalf, -\vecbeta) >0$ (recall the explicit dependence of $C$ on $\theta_k$ given in Remark \ref{remark-size-constant} and the definition of $m_k$).  In other words, the constant $M$ does not depend on the group element $\gamma$ we used to  translate the point $(x+iN^{-2},0;\sve{\vecalf}{-\vecbeta})$ back to the fundamental domain $\calF_\Gamma$.

    By Lemma \ref{kappa_shift} we have that \begin{align}\label{bound-of-kappa_eta(f_N)}
        \kappa_\eta (f_N) \leq 3^\eta(1 + \tfrac{1}{N^2}\|\vecbeta\|^2)^{\eta/2}\kappa_\eta (f).
    \end{align}
        Finally, combining \eqref{4-20},  \eqref{bound-with-uniform-constant-over-orbit}--\eqref{bound-of-kappa_eta(f_N)}, and using the trivial lower bound $N\geq1$, we obtain
    \begin{align}
        |S_N^f(x;\vecalf,\vecbeta)|  \leq M\, 3^\eta  \,  (1 + \|\vecbeta\|^2)^{\eta/2}\,\kappa_\eta (f) N^{k/2}.\label{regular-f-bound-inf_eta}
    \end{align}
    Note that if $f\in \Si_{\eta}(\R^k)$ for some $\eta$ then $f$ also belongs to any $\Si_{\eta'}(\R^k)$ with $\eta'\leq \eta$. Therefore \eqref{regular-f-bound-inf_eta} implies \eqref{bounded-sums-chapter-5-statement} with 
    \begin{align}\label{implied-constant}
        \inf_{\eta'\in(k,\eta]}\lp M(k,m_k,\eta')\,3^{\eta'} \,(1 + \|\vecbeta\|^2)^{\eta'/2}\,\kappa_{\eta'} (f) \rp
    \end{align}
    as implied constant.
\end{proof}
\begin{remark}\label{remark:remove-explicit-dependence-upon-vecbeta}
    If we add the minor assumption  that $N \geq \|\vecbeta\|$ in Theorem \ref{bounded-sums-chapter-5}, then we may remove the explicit $\vecbeta$-dependence of the implied constant in \eqref{bounded-sums-chapter-5-statement}. In fact, bounding $1 + \tfrac{1}{N^2}\|\vecbeta\|^2$ above by $2$ in \eqref{bound-of-kappa_eta(f_N)} instead of using the trivial bound $N\geq1$, we may replace \eqref{implied-constant} with \begin{align}\inf_{\eta'\in(k,\eta]} M(k,m_k,\eta')\,3^{\eta'}2^{\eta'/2}\,\kappa_{\eta'} (f).
    \end{align}Naturally, the implied constant still depends implicitly on $(\vecalf,\vecbeta)$ via $m_k(\vecalf,-\vecbeta)$. 
\end{remark}

%% file: study-of-the-orbits.tex
\subsection{Study of  $\Gamma$-Orbits  
and Conclusions}\label{study-of-orbits}
Set $\vecxi_1 := (\xi_{11}, \cdots, \xi_{1n})^{\top}$ and $\vecxi_2 := (\xi_{21}, \cdots, \xi_{2n})^{\top}$. Define
\begin{align}
    \phi&: \R^{2k} \to \R^{2k}\\
    \phi&: \sve{\vecxi_1}{\vecxi_2} \mapsto \begin{pmatrix} \sve{\xi_{11}}{\xi_{21}} \\ \vdots \\ \sve{\xi_{1n}}{\xi_{2n}}\end{pmatrix}.
\end{align}
Then
\begin{align}
    \phi \lp\smatr{a}{b}{c}{d} \cdot \sve{\vecxi_1}{\vecxi_2}\rp = \begin{pmatrix} M\sve{\xi_{11}}{\xi_{21}} \\ \vdots \\ M\sve{\xi_{1n}}{\xi_{2n}}\end{pmatrix}.
\end{align}
Therefore $\Orb_{\Gamma}\!\sve{\vecxi_1}{\vecxi_2}$ is determined by the orbits of $\sve{\xi_{1i}}{\xi_{2i}}$ under the action of $\Gamma_\theta$ on $\R^2/ \Z^2$, $\Orb_{\Gamma_{\theta}}\!\sve{\xi_{1i}}{\xi_{2i}}$, $1 \leq i \leq n$. Recalling the definition \eqref{m-function}, we can use the previous observation to can guarantee that $m_k(\vecxi_1, \vecxi_2) > 0$  by checking that  $m_1(\xi_{1i}, \xi_{2i})$ is positive for at least one index $1 \leq i \leq n$. 
To this end, let $L_U := \{\sve{\xi_1}{\xi_2} \in \R^2: \xi_2 = 0\}$, and $L^\pm_V := \{\sve{\xi_1}{\xi_2} \in \R^2 : \xi_2 = \xi_1 \pm \tha \}$. For $(\alpha, \beta) \in \R^2$ we define
\begin{align}
   U\toab &:= L_U \cap \Orb_{\Gamma_\theta}\sve{\alpha}{\beta}\\
   V\toab &:= L_V \cap \Orb_{\Gamma_\theta}\sve{\alpha}{\beta}.
\end{align}
We have the following 
\begin{proposition}
    \label{set-K-classification}
    Let $(\alpha, \beta) \in \Q^2$ where $\{\alpha\} = \tfrac{a}{q}$ and $\{\beta\} = \tfrac{b}{q}$ with $\gcd (a,b,q) = 1$ and $q > 1$. Then $|U\toab| = |V\toab| = 0$ if and only if $a$ and $b$ are both odd, and $q = 2m$ where $m$ is odd.
\end{proposition}
\begin{proof}
    This follows immediately from Propositions 3.2.3, 6.2.1 and 6.1.1 in \cite{Cellarosi-Osman-rational-tails}. 
\end{proof}
Note that $|\Orb_{\Gamma_\theta}\!\sve{\alpha}{\beta}| < \infty$ for any $(\alpha,\beta) \in \Q^2$. It is important to observe that for $(\alpha,\beta) \in \Q^2$ we have that $m_1(\alpha,\beta)>0$ if and only if the $\Gamma_\theta$-orbit of $\sve{\alpha}{\beta}$ does not intersect neither $L_U$ nor $L_V^{\pm}$. Hence, Proposition \ref{set-K-classification} gives a classification of all rational pairs $(\alpha, \beta) \in \Q^2$ for which $m_1(\alpha, \beta) >0$.

The orbit $\Orb_{\Gamma_\theta} \!\sve{\alpha}{\beta}$ for $(\alpha,\beta) \in \Q^2$ also exhibits various symmetries: 
\begin{proposition}\label{orbit-symmetries}
For any $(\alpha,\beta) \in \Q^2$, we have that
\begin{align}
    m_1(\alpha,\beta) = m_1(\alpha, -\beta) = m_1(-\alpha,\beta) = m_1(-\alpha,-\beta).
\end{align}
\end{proposition}
\begin{proof}
    This is immediate from Lemma 3.3.1 in \cite{Cellarosi-Osman-rational-tails}.
\end{proof}

We can now state and prove our Main Theorem, which generalises Theorem  \ref{thm:main-thm-intro} from Schwartz cut-offs to to regular cut-offs. It is basically a version of Theorem \ref{bounded-sums-chapter-5}, where the assumption that $m_k(\vecalf,-\vecbeta)>0$ is replaced by a concrete, easy to check, sufficient condition (coming from Proposition \ref{set-K-classification}).


\begin{theorem}[Main Theorem, for regular cut-offs]\label{thm:main-theorem}
    Let $f \in \Si_\eta(\R^k)$ with $\eta>k$. Let $\vecalf, \vecbeta \in \R^k$, where $\vecalf = (\alpha_1, \cdots, \alpha_k)$ and $\vecbeta = (\beta_1, \cdots, \beta_k)$. Suppose that at least one of $(\alpha_i,\beta_i)$ for $1\leq i\leq k$ can be written as $(\alpha_i,\beta_i)= (\tfrac{a}{2m}, \tfrac{b}{2m}) \in \Q^2$ with  $\gcd(a,b, m) = 1$ and $a$, $b$, and $m$ all odd. Let $\beta$ the smallest (in absolute value) of the $\beta_i$'s corresponding to such pairs. Then 
    \begin{align}\label{statement-main-theorem-regular}
        |S_N^f (x;\vecalf, \vecbeta)| \ll_{k,m,\eta,\beta,f}
        N^{k/2}.
    \end{align} 
\end{theorem}

\begin{proof}
    Suppose that, without loss of generality, $(\alpha_1,\beta_1) = (\tfrac{a}{2m}, \tfrac{b}{2m}) \in \Q^2$ is such that $\gcd(a,b, m) = 1$, and $a$, $b$, and $m$ are all odd. By Propositions \ref{set-K-classification} and \ref{orbit-symmetries} we are guaranteed that $m_1(\alpha_1, -\beta_1) = m_1(\alpha_1,\beta_1) > 0 $. Clearly $m_k(\vecalf,-\vecbeta) \geq m_1(\alpha_1,\beta_1) > 0$. It is easy to check that $m_1(\alpha_1,\beta_1)\geq\frac{1}{2m}$ (see Lemma 5.0.2 in \cite{Cellarosi-Osman-rational-tails}) and so by Theorem \ref{bounded-sums-chapter-5} we have the result.
\end{proof}
\begin{proof}[Proof of Theorem \ref{thm:main-thm-intro}]
If $f$ is of Schwartz class, then we can apply Theorem \ref{thm:main-theorem} for every $\eta>k$. Equivalently, in the case of Scwhartz cut-offs, 
the dependence upon $\eta$ in
\eqref{implied-constant} ---and hence in \eqref{statement-main-theorem-regular}---
is removed by taking the infimum over all $\eta'>k$.
\end{proof}

\begin{remark}\label{remark:k=1}
    If $k = 1$ then
    \begin{align*}
        C(k, \xi_2, \eta) &= 2^{\tfrac{2\eta - 1}{4}} \lp \frac{1}{\theta_1(\xi_2)^\eta} + \sum_{n = 1}^{\infty} \lp \frac{1}{(n - \theta_1(\xi_2))^\eta} + \frac{1}{(n + \theta_1(\xi_2))^\eta}\rp\rp,
    \end{align*}
    and similarly for $C(k, \hat\xi_2, \eta$). Therefore, we may bound $M(k,m_1,\eta')$ in terms of the Riemann zeta function evaluated at $\eta'$ and in terms of the  denominator $2m$ used to uniquely write $(\alpha,\beta)$ as $\alpha=\frac{a}{2m}$, $\beta=\frac{b}{2m}$ with 
    $\gcd(a,b,m)=1$. 
 We obtain
    \begin{align}
         \inf_{\eta'\in(1,\eta]} \lp  2^{\frac{2\eta' - 1}{4}}3^{\eta'}(2m + 2\zeta(\eta')(2^{\eta'} - 1))
        (1 + \beta^2)^{\eta'/2}\kappa_{\eta'} (f) \rp\label{eq:k=1_constant}
    \end{align}
    as an implied constant for \eqref{statement-main-theorem-regular}. For Schwartz cut-offs, such as those in Corollary \ref{cor:main-cor-intro}, the infimum in \eqref{eq:k=1_constant} is taken over all $\eta'>1$. See Section \ref{jacobi-theta-bound-section}.

\end{remark}

\begin{remark}\label{remark:additional-assumption2}
    Similarly to Remark \ref{remark:remove-explicit-dependence-upon-vecbeta}, if we make the additional mild assumption that $N\geq\|\vecbeta\|$ in Theorem \ref{thm:main-theorem}, then we can remove the explicit dependence upon $\vecbeta$ in \eqref{statement-main-theorem-regular}. In the case where $k=1$ the factor $(1+\beta^2)^{\eta'/2}$ in \eqref{eq:k=1_constant} may then be replaced by $2^{\eta'/2}$. This is done in Theorem \ref{thm:jacobi-theta-bound-intro} since $(\alpha,\beta)$ are assumed to be in the torus $\in\R^2/\Z^2$.  
\end{remark}

%% file: theta-function-bounds.tex
\section{Applications}\label{section:applications}

\subsection{Uniform Bounds for $\vtheta_3$}\label{jacobi-theta-bound-section}
    The classical Jacobi theta function $\vtheta_3(z,w)$ is defined in \eqref{def-classical-jacobi-theta3}
    for $z, w \in \C$ with $\Im(w) > 0$. It is often written as $q$-series (where $q = e^{-\pi i w}$) as 
    \begin{align}
      \vtheta_3 (z; q) = \sum_{n \in \Z} e(nz)q^{n^2} = 1 + 2\sum_{n = 0}^{\infty}\cos (2nz)q^{n^2},
    \end{align}
    where $|q| < 1$. Note that if $\Im(w)\to0^+$, then $|q|\to1^-$.
    Considering the Schwartz function $f(u) = e^{-\pi u^2}$ in \eqref{key-relation-theta-sum-to-function} we see that 
    \begin{align}
      |\Theta_f(x + iy, 0; \sve{\xi_1}{0})| &= y^{1/4} \labs\sum_{n \in \Z} e^{2 \pi i \xi_1 n + \pi i \lp x + iy\rp n^2}\rabs\\
      &=  y^{1/4} |\vartheta_3 (\xi_1, x + iy)|.
    \end{align}
    Taking $y = \varepsilon = N^{-2}$ and $\xi_1 = \alpha + \beta x$ gives
    \begin{align}
      \frac{1}{\sqrt N}|S_N^f(x;\alpha,\beta)|= \varepsilon^{\frac{1}{4}} \labs\vtheta_{3}\!\lp \alpha + \beta x, x + i \varepsilon 
      \rp\rabs.
    \end{align}


We are now ready to prove Theorem \ref{thm:jacobi-theta-bound-intro}
where the implied constant in \eqref{statement-thm-bound-theta3} can be taken as \begin{align}\label{statement-thm-bound-theta3-implied-constant}
\min_{1<\eta<4/3} \lp 2^{-\frac{1}{4}}\, 6^\eta  \,(2m + 2\zeta(\eta)(2^{\eta} - 1))\rp.
\end{align}
\begin{proof}[Proof of Theorem \ref{thm:jacobi-theta-bound-intro}]
Recall \eqref{rotation-action-on-L2}-\eqref{phi-transform-definition}.
We claim that for  $f(u)  = e^{-\pi u^2}$ and every $\phi\in\R$ we have \begin{align}\label{phi-transform-gaussian}
|f_{\phi}(w)| = e^{-\pi w^2}.
\end{align}
This is obvious if $\phi\equiv0\:(\bmod\:\pi)$. For all other $\phi$'s, we first use the  identity (see \cite{Abramowitz-Stegun} 7.4.32)
    \begin{align}
        \int_{-\infty}^{\infty}e^{A + Bt + Ct^2}\:\de t = \frac{e^{A - \frac{B^2}{4C}}\sqrt \pi}{\sqrt{-C}}
    \end{align}
    for $A,B,C \in \C$ with $\Re(C) < 0$.  We have 
    \begin{align}
        |f_\phi (w)| &= \labs\frac{1}{\sqrt{|\sin \phi|}}\int_{-\infty}^{\infty}e\!\lp\frac{\tha (w^2 + t^2)\cos \phi + wt}{\sin \phi}\rp e^{-\pi t^2} \: \de t \rabs\\
        &= \labs\frac{\sqrt \pi}{\sqrt{|\sin \phi|}} \frac{e^{\pi i w^2 \cot \phi - \frac{-2\pi i w \csc \phi}{4(-\pi + i\pi\cot \phi)}}}{\sqrt{-(-\pi + i\pi \cot \phi)}}\rabs
    = \labs \frac{1}{\sqrt{|\sin \phi|}}\frac{e^{-\pi w^2}}{\sqrt{1 - i \cot \phi}}\rabs.
    \end{align}
    Note that $|(1 - i \cot \phi)^{-1/2}| = |(1 - i \cot \phi)(1 + i\cot \phi)|^{-1/4} = |\csc^2\phi|^{-1/4}= |\csc \phi|^{-1/2}$ and so
    \begin{align}
        \labs\frac{1}{\sqrt{|\sin \phi|}\sqrt{1 - i\cot \phi}}\rabs = 1,
    \end{align}
    yielding \eqref{phi-transform-gaussian}.
    Now recall  and \eqref{kappa-norm-def}  let  $\eta > 1$. We claim that
    \begin{align}\label{formula-kappa-eta-gaussian}
        \kappa_{\eta}(f) = \begin{dcases*} 1 & if $1 < \eta < 2 \pi$ \\ e^{\pi}(\tfrac{\eta}{2 \pi})^{\eta/2}e^{-\eta/2} & if $\eta \geq 2 \pi$.\end{dcases*}
    \end{align} 
    By \eqref{phi-transform-gaussian}
   we have
    \begin{align}
    \kappa_{\eta}(f) = \sup_{\phi, w}(1 + w^2)^{\eta/2} |f_{\phi}(w)|= \max_{w}(1 + w^2)^{\eta/2}e^{-\pi w^2}.
    \end{align}
    By finding the zeroes of the derivative $\tfrac{\de}{\de w} (1 + |w|^2)^{\eta/2}e^{-\pi w^2} = e^{-\pi w^2}(w (w^2 + 1)^{\tfrac{\eta}{2} - 1}(\eta - 2\pi (w^2 + 1)))$, it is not hard to see that if $1 < \eta \leq 2 \pi$, then $\max_{w}(1 + w^2)^{\eta/2}e^{-\pi w^2}$ occurs at $w = 0$ and equals 1. 
    If $\eta > 2 \pi$, then  $\max_{w}(1 + w^2)^{\eta/2}e^{-\pi w^2}$ occurs at $w = \pm \sqrt{\tfrac{\eta}{2 \pi} - 1}$ and equals 
   $e^{\pi}(\tfrac{\eta}{2})^{\eta /2}e^{-\eta /2}$. 

   Now let us apply Corollary \ref{cor:main-cor-intro} with $f(u)=e^{-\pi u^2}$ and $N=\varepsilon^{-\frac{1}{2}}$. The implied constant in  \eqref{cor:main-cor-intro-bound} can be taken as
   $\inf_{\eta>1}h(m,\eta)$, where $h(m,\eta):=3^\eta 2^{\eta-\frac{1}{4}}(2m + 2\zeta(\eta)(2^{\eta} - 1))\kappa_{\eta} (f)$
   (see Remarks \ref{remark:k=1}-\ref{remark:additional-assumption2}).
   Using \eqref{formula-kappa-eta-gaussian}
and recalling that the $\zeta$-function $\eta\mapsto\zeta(\eta)$  has a simple pole at $\eta=1$,  we see that $\lim_{\eta\to1^+}  h(m,\eta)=\lim_{\eta\to\infty}h(m,\eta)=\infty$ for every $m\geq1$ and hence the infimum is a minimum, attained at some $\eta^*(m)>1$. We leave it to the reader to verify that such point of minimum is unique. Since $m\mapsto h(m,\eta)$ is increasing, we have $\eta^*(m)\leq \eta^*(1)$. It is clear that $\eta^*(1)\leq 2\pi$ and we can numerically estimate $\eta^*(1)\approx1.263<\frac{4}{3}$, i.e. when taking the minimum we can restrict $\eta$ to the interval $(1,\frac{4}{3})$,  where $\kappa_\eta(f)=1$, thus obtaining \eqref{statement-thm-bound-theta3-implied-constant}. 
\end{proof}

\begin{figure}[h!]
    \centering
    \includegraphics[width=16.25cm]{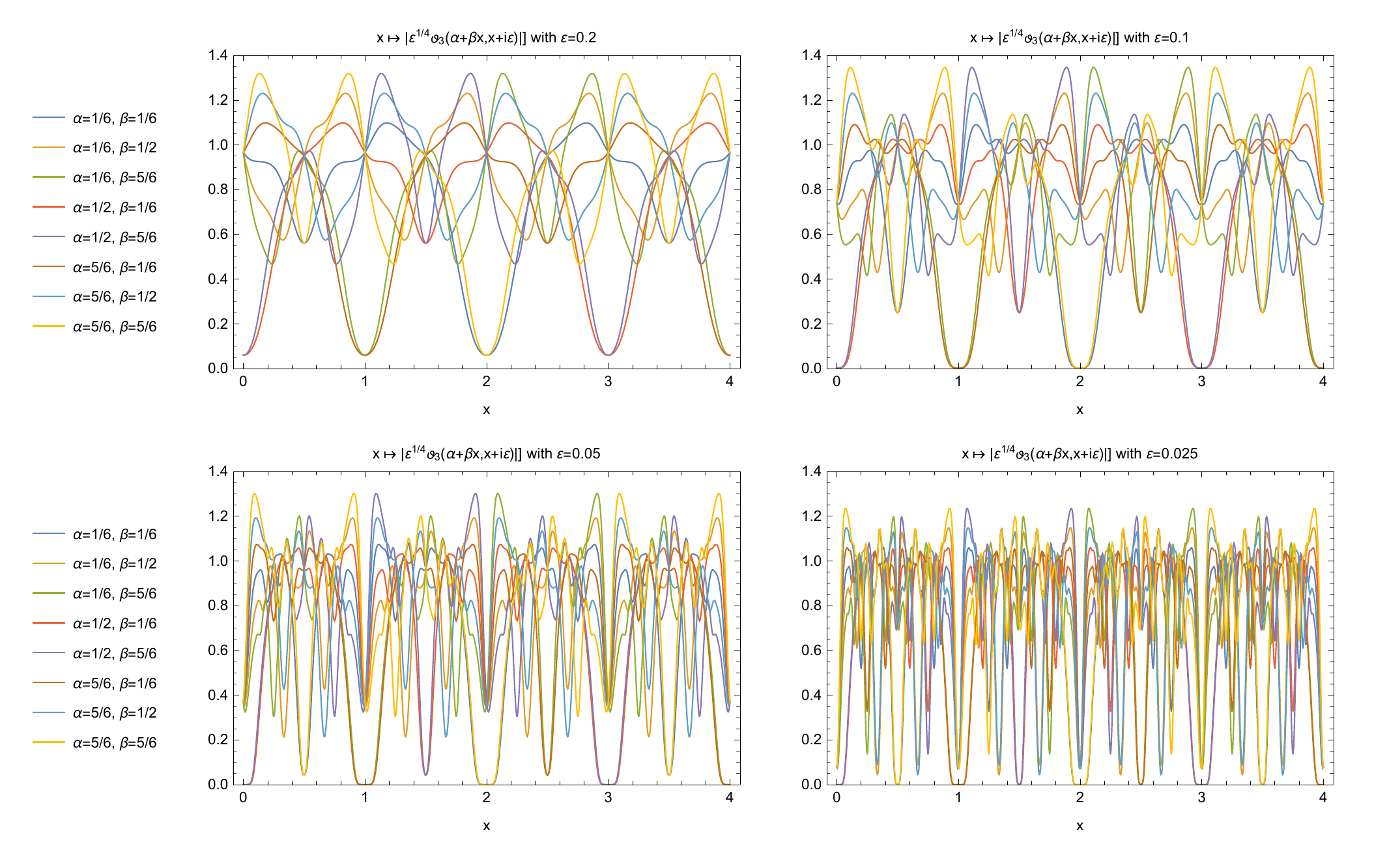}
    \caption{Illustration of Theorem  \ref{thm:jacobi-theta-bound-intro} for $m=3$. There are $J_2(3)=8$ rational pairs $(\alpha,\beta)\in\R^2/\Z^2$ satisfying $\alpha=\frac{a}{2m}$, $\beta=\frac{b}{2m}$ with $\gcd(a,b,m)=1$ and $a,b,m$ all odd, see Remark \ref{abundance-of-pairs-with-same-bound}. The corresponding 8 functions $[0,4]\ni x\mapsto \varepsilon^{\frac{1}{4}} \labs\vtheta_{3}\!\lp \alpha + \beta x, x + i \varepsilon \rp\rabs\in\R_{\geq0}$ are shown for various values of $\varepsilon\in\{0.2, 0.1,0.05,0.025\}$. For one of these pairs, $(\frac{1}{2},\frac{1}{6})$, the function $x\mapsto \varepsilon^{\frac{1}{4}} \labs\vtheta_{3}\!\lp \alpha + \beta x, x + i \varepsilon \rp\rabs\in\C$ is shown in Figure \ref{fig:six-figures-theta3-example}.}
    \label{fig:many-theta3-same-bound}
\end{figure}

\begin{remark}\label{abundance-of-pairs-with-same-bound}
    Note the bound \eqref{statement-thm-bound-theta3} is the same for all rational pairs $(\frac{a}{2m},\frac{b}{2m})$ with $\gcd(a,b,m)=1$ and $a, b, m$ all odd, since \eqref{statement-thm-bound-theta3-implied-constant} only involves the denominator $2m$. It is easy to see that, given $m$ odd, the number of such pairs in $\R^2/\Z^2$ is between $\frac{8m^2}{\pi^2}$ and $m^2$, their exact number being given by the second Jordan totient function $J_2(m)=m^2\prod_{p|m}\lp1-p^{-2}\rp$, see Remark 7.1.2 in \cite{Cellarosi-Osman-rational-tails}. Figure \ref{fig:many-theta3-same-bound} illustrates the uniformity of the bound \eqref{statement-thm-bound-theta3} across all rational pairs in $\R^2/\Z^2$ with denominator $6$ that satisfy the hypotheses of Theorem \ref{thm:jacobi-theta-bound-intro}.
\end{remark}

\subsection{Estimates for $S_N$ independent of $x$}
Recall that 
    $S_N(x; \alpha,\beta) = \sum_{n=1}^N e((\tha n^2 + \beta n)x + \alpha n)$.
For fixed $N$ we can find a regular cut-off function $\chi^{(N)} \in S_\eta(\R^k)$ with $\eta=2$,  supported on $[0,1 + \tfrac{1}{N}]$, and $\chi^{(N)}(w) = 1$ for $w \in [\tfrac{1}{N}, 1]$. An explicit construction is given in \cite{Cellarosi-Marklof}, equation (3.9). In particular,
\begin{align}
    \tfrac{1}{\sqrt N}|S_N(x; \alpha,\beta)| = |\Theta_{\chi^{(N)}}(x + \tfrac{1}{N^2}, 0; \sve{\alpha + \beta x}{0})|. \label{eq:S_N-and-Theta}
\end{align}
Proposition 6.5 in \cite{Cellarosi-Griffin-Osman-error-term} shows that 
\begin{align}\label{kappa_eta_chi_N}
    \kappa_\eta(\chi^{(N)}) \ll N^{\eta - 1}.
\end{align}
We obtain the following corollary of Main Theorem \ref{thm:main-theorem}.

\begin{cor}\label{cor:S_N-bound}
    Suppose $(\alpha, \beta) = (\tfrac{a}{2m}, \tfrac{b}{2m}) \in \Q^2$ where $\gcd(a,b,m) = 1$ and $a,b,$ and $m$ are all odd. Then
    \begin{align}
        |S_N(x; \alpha,\beta)|\ll_{\varepsilon} 
        N^{\tha + \varepsilon}
        \label{eq:S_N-bound}
    \end{align}
    for every $0 < \varepsilon \leq 1$ and every $x, N$.
\end{cor}

\begin{proof}
Using \eqref{eq:S_N-and-Theta},  Theorem \ref{thm:main-theorem} with implied constant \eqref{eq:k=1_constant}, and \eqref{kappa_eta_chi_N}, we obtain
    \begin{align}
        |S_N(x; \alpha,\beta)| &\leq \lp   2^{\frac{2\eta - 1}{4}}3^\eta(2q + 2\zeta(\eta)(2^{\eta} - 1))(1 + \beta)^{\eta}\kappa_\eta (\chi^{(N)}) \rp N^{\tha} \\
        &\ll \lp   2^{\frac{2\eta - 1}{4}}3^\eta(2q + 2\zeta(\eta)(2^{\eta} - 1))(1 + \beta)^{\eta} \rp N^{\tha + (\eta - 1)}.
    \end{align}
    Taking $\varepsilon=\eta-1$ close to $0^+$ gives the result.
\end{proof}

\begin{remark}
    Notice that as $\varepsilon \to 0$, the implied constant in \eqref{eq:S_N-bound} goes to $\infty$ due to the presence of the Riemann $\zeta$-function in the implied constant.
\end{remark}

\begin{remark}\label{Remark-FF-bound}
    As mentioned in the introduction, this estimate is not optimal. For instance, it follows from Corollary 1.2 in \cite{Flaminio-Forni2006} (cf. also \cite{Cosentino-Flaminio},\cite{Marklof-Welsh-higher-rank-theta-sums}) that there is a full measure set $\calA \subseteq [0,1]$ such that for any $x \in \calA$
    \begin{align}\label{bound-flaminio-forni}
        |S_N(x; \alpha,\beta)| \ll _{x,\varepsilon}\sqrt N (\log N)^{\frac{1}{4} + \varepsilon}.
    \end{align}
    for every $N$. For typical $x$, this bound is certainly better than \eqref{eq:S_N-bound}. On the other hand, for the class of rational pairs $(\alpha,\beta)$ we consider, the constant implied in \ref{bound-flaminio-forni} does depend on $x$, while the one in \eqref{eq:S_N-bound} is uniform in $x$.
\end{remark}

\section*{Acknowledgements} The first author acknowledges the support from the NSERC Discovery Grant 
RGPIN-2022-04330. The results presented here were partly developed in the PhD thesis of the second author. We wish to thank Jens Marklof, Ram M. Murty, and Brad Rodgers for several fruitful discussions on the subject of this work.